\documentclass[a4paper,10pt]{article}
\usepackage{amssymb,amsmath,amsthm}
\usepackage{fullpage}
\usepackage{hyperref}
\usepackage{eucal}
\usepackage{graphicx}
\newcommand{\ie}{\emph{i.e.}}

\newcommand{\cf}{\emph{cf.}}
\newcommand{\Real}{\mathbb{R}}
\newcommand{\Com}{\mathbb{C}}

\newcommand{\supp}{\mathop{\mathrm{supp}}\nolimits}
\newcommand{\Dom}{\mathsf{D}}

\newcommand{\tr}{\mathop{\mathrm{tr}}\nolimits}
\newcommand{\eps}{\varepsilon}
\newcommand{\sii}{L^2}
\newcommand{\Hilbert}{\mathcal{H}}
\newcommand{\der}{\mathrm{d}}
\newtheorem{Theorem}{Theorem}
\newtheorem{Lemma}{Lemma}
\newtheorem{Proposition}{Proposition}
\newtheorem{Corollary}{Corollary}

\theoremstyle{definition}
\newtheorem{Remark}{Remark}

\usepackage[normalem]{ulem}
\usepackage{color}
\definecolor{DarkBlue}{rgb}{0,0.1,0.7}

\newcommand\soutD{\bgroup\markoverwith
{\textcolor{DarkBlue}{\rule[.01ex]{2pt}{1pt}}}\ULon}
\newcommand{\Hm}[1]{\leavevmode{\marginpar{\tiny%
$\hbox to 0mm{\hspace*{-0.5mm}$\leftarrow$\hss}%
\vcenter{\vrule depth 0.1mm height 0.1mm width \the\marginparwidth}%
\hbox to
0mm{\hss$\rightarrow$\hspace*{-0.5mm}}$\\\relax\raggedright #1}}}

\begin{document}
%
%-------%
% TITLE %
%-------%
%------------------------------------------%
%------------------------------------------%
\title{\textbf{\LARGE
Location of eigenvalues of three-dimensional non-self-adjoint Dirac operators
}}
\author{Luca Fanelli\,$^a$ \ and \ David Krej\v{c}i\v{r}{\'\i}k\,$^b$}
\date{\small 
\begin{quote}
\emph{
\begin{itemize}
\item[$a)$] 
Dipartimento di Matematica, SAPIENZA Universit\`a di Roma,
P.~le Aldo Moro 5, 00185 Roma;
fanelli@mat.uniroma1.it.%
\item[$b)$] 
Department of Mathematics, Faculty of Nuclear Sciences and 
Physical Engineering, Czech Technical University in Prague, 
Trojanova 13, 12000 Prague 2, Czechia;
david.krejcirik@fjfi.cvut.cz.%
\end{itemize}
}
\end{quote}
19 September 2018}
\maketitle
\begin{abstract}
\noindent
We prove the absence of eigenvaues of the three-dimensional Dirac operator 
with non-Hermitian potentials in unbounded regions of the complex plane
under smallness conditions on the potentials in Lebesgue spaces.
Our sufficient conditions are quantitative and easily checkable.
\bigskip
\begin{itemize}
\item[\textbf{Keywords:}]
Dirac operator, complex potential, non-self-adjoint perturbation,
pseudo-Friedrichs extension, Birman-Schwinger principle,
absence of eigenvalues. 
\item[\textbf{MSC (2010):}]
Primary: 35P15, 35J99, 47A10, 47F05, 81Q12.
\end{itemize}
\end{abstract}
%
%------------------------------------------%
%------------------------------------------%
 
%---------------------%
\section{Introduction}
%---------------------%
%
Let us consider a relativistic quantum particle of spin $\frac{1}{2}$
and mass $m \geq 0$ in~$\Real^3$,
subject to an external electric field described by a potential~$V$.
The dynamics is governed by the Dirac Hamiltonian
\begin{equation}\label{Dirac}
  H_{V} := - i \;\! \alpha \cdot \nabla  + m \, \alpha_4 + V
\end{equation}
acting in the Hilbert space of spinors 
$\Hilbert := \sii(\Real^3;\Com^4)$.
Here $\alpha := (\alpha_1,\alpha_2,\alpha_3)$
with~$\alpha_\mu$ being the usual $4 \times 4$ Hermitian Dirac matrices 
satisfying the anticommutation rules
\begin{equation}\label{anti}
  \alpha_\mu \alpha_\nu + \alpha_\nu \alpha_\mu
  = 2 \delta_{\mu\nu} I_{\Com^4} 
\end{equation}
for $\mu,\nu\in\{1,\dots,4\}$ 
and the dot denotes the scalar product in~$\Real^3$.
Motivated by a growing interest in non-self-adjoint operators
in quantum mechanics (\cf~\cite{Bagarello-book} for a mathematical overview),
we proceed in a greater generality by allowing~$V$
to be a possibly non-Hermitian $4 \times 4$ matrix in~\eqref{Dirac}.

In the traditional self-adjoint case 
(\ie~$V$ is a real scalar multiple of the diagonal matrix),
the literature on spectral properties of~$H_V$ is enormous
and we limit ourselves to quoting the classical Thaller's monograph~\cite{Thaller}.
For~$V$ being matrix-valued and possibly non-Hermitian,
a systematic study of the spectrum of~$H_V$ was initiated
by the pioneering work of Cuenin, Laptev and Tretter 
\cite{Cuenin-Laptev-Tretter_2014} in the one-dimensional setting
and followed by \cite{Cuenin_2014,Enblom_2018,Cuenin-Siegl_2018}.
Some spectral aspects in the present three-dimensional 
situation are covered by \cite{Dubuisson_2014,Sambou_2016,Cuenin_2017}.

In the field-free case (\ie\ $V=0$),
it is well known that $\sigma(H_{0}) = (-\infty,-m] \cup [+m,+\infty)$
and that the spectrum is purely continuous.
That is, the residual spectrum is empty and there are no eigenvalues.
The objective of this paper is to derive quantitative smallness 
conditions on the potential~$V$,
which guarantee that the spectrum of~$H_{V}$ remains purely continuous,
at least in certain regions of the complex plane.
Denoting by $|V(x)|$ the operator norm of the matrix $V(x) \in \Com^{4 \times 4}$
for almost every fixed $x \in \Real^3$, 
the smallness is measured through Lebesgue norms of 
the real-valued function~$|V|$.

We present two types of results in this paper.
The first reads as follows:

\begin{Theorem}\label{Thm1}
Assume $|V| \in L^3(\Real^3)$.
If  
\begin{equation}\label{hypothesis1}
  C \,f(\lambda,m) \,\||V|\|_{L^3(\Real^3)} 
  < 1
\end{equation}
with
$$
  C := \left(\frac{\pi}{2}\right)^{1/3} \sqrt{1+e^{-1}+2e^{-2}}
  \qquad \mbox{and} \qquad
  f(\lambda,m) := \sqrt{ 
  1 + \frac{(\Re\lambda)^2}{\big(\Re\sqrt{m^2-\lambda^2}\big)^2} 
  }
  \,,
$$
then $\lambda \not\in \sigma_\mathrm{p}(H_V)$.
\end{Theorem}

\begin{figure}[h]
\begin{center}
\includegraphics[width=0.5\textwidth]{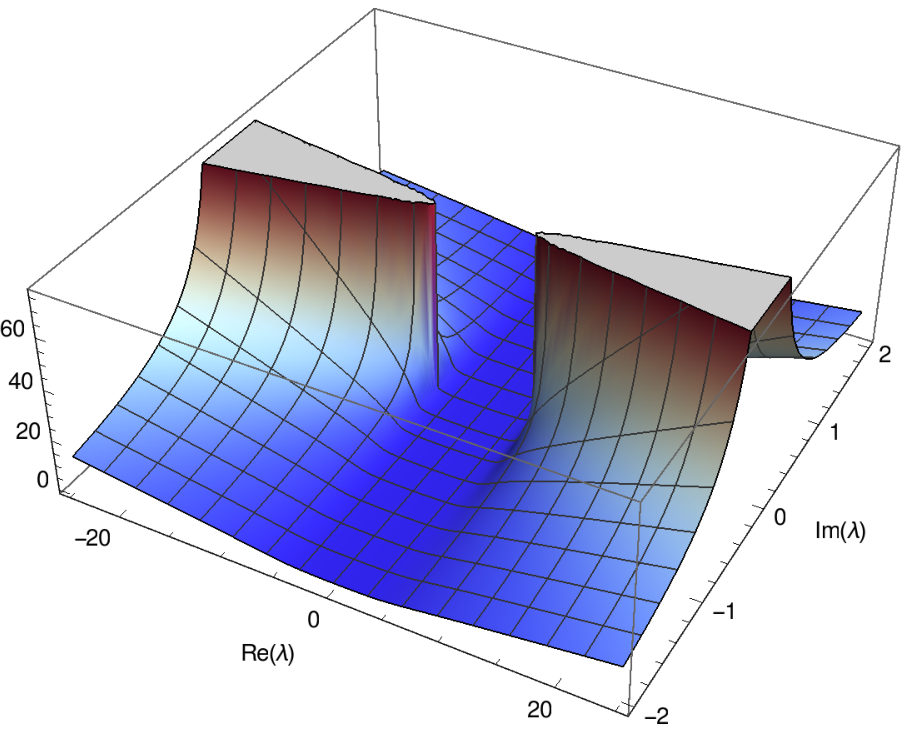}
\
\includegraphics[width=0.4\textwidth]{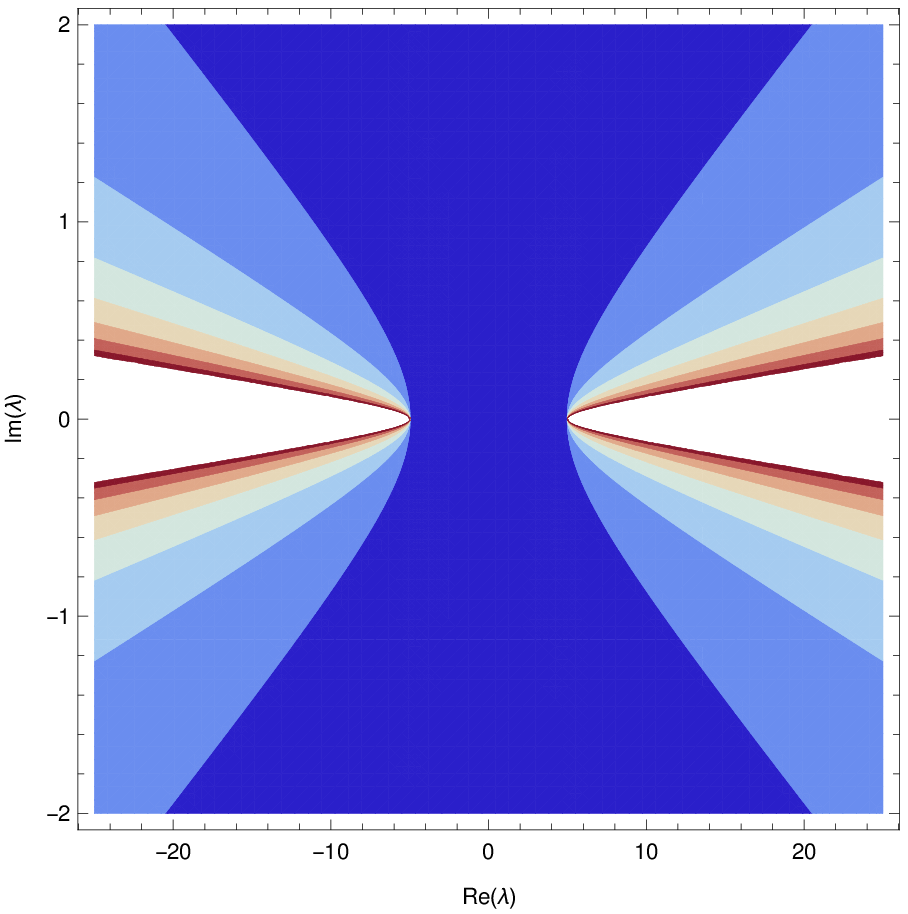}
\includegraphics[width=0.5cm]{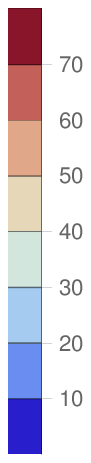}%
\caption{Graph of the surface $\lambda \mapsto f(\lambda,m)$ 
and its contour plot for $m=5$.}\label{Fig}
\end{center}
\end{figure}

The hypothesis~\eqref{hypothesis1} is essentially a smallness requirement about~$V$.
In view of Kato's general smoothness theory \cite{Kato_1966}
(see also \cite[Sec.~XIII.7]{RS4}), it is not surprising that a result
of this type should hold for weakly coupled potentials.
The strength of our result lies in that the condition is explicit
and easy to check in applications.
Notice that $C \approx 1.5$ and that $f(\lambda,m)$ is finite 
if, and only if, $\lambda \not\in (-\infty,-m]\cup[+m,+\infty)$.
Consequently, the condition~\eqref{hypothesis1} is obeyed for such~$\lambda$
whenever the norm $\||V|\|_{L^3(\Real^3)}$ is sufficiently small.  
Let us also remark that 
$f(\lambda,m) \sim |\lambda|/|\Im(\lambda)|$ as $|\Im(\lambda)| \to +\infty$
and in fact $f(\lambda,0)=|\lambda|/|\Im(\lambda)|$ in the massless case. 
We illustrate the dependence of $f(\lambda,m)$ on~$\lambda$ in Figure~\ref{Fig}.

\begin{Remark}
An analogous result holds (with an unspecified constant)
provided that $L^{3}(\Real^3)$ is replaced by 
the Lorentz space $L^{3,\infty}(\Real^3)$.
Since this space in particular contains the Coulomb potential
$V_Z(x) := -Z/|x|$, which creates discrete eigenvalues 
in the gap $(-m,+m)$ whenever $Z>0$, \cf~\cite[Sec.~7.4]{Thaller}, 
we see that the presence of a $\lambda$-dependent function $f(\lambda,m)$
diverging as $\lambda \to \pm m$ is in fact unavoidable in~\eqref{hypothesis1}.
\end{Remark}

To get a more uniform control over the spectrum 
and exclude possible eigenvalues embedded in 
the essential spectrum $(-\infty,-m]\cup[+m,+\infty)$,
we have to strengthen the condition about the potential~$V$. 

\begin{Theorem}\label{Thm2}
Assume $|V| \in L^3(\Real^3) \cap L^{3/2}(\Real^3)$.
If  
\begin{equation}\label{hypothesis2}
  C \, \||V|\|_{L^3(\Real^3)} 
  + C' \, |\Re \lambda| \,
  \||V|\|_{L^{3/2}(\Real^{3})} 
  < 1
  \,,
\end{equation}
where~$C$ is the same constant as in Theorem~\ref{Thm1}
and 
$$
  C' := \frac{2^{17/6}}{3\pi^{2/3}}
  \,,
$$
then $\lambda \not\in \sigma_\mathrm{p}(H_V)$.
\end{Theorem}

In this case, 
given a potential~$V$ with sufficiently small norm $\||V|\|_{L^3(\Real^3)}$,
the hypothesis~\eqref{hypothesis2} excludes the existence of eigenvalues 
in thin tubular neighbourhoods of the imaginary axis,
with the thinness determined by the norm $\||V|\|_{L^{3/2}(\Real^3)}$.
Notice that $C' \approx 1.1$ 
and that the condition~\eqref{hypothesis2} is $m$-independent.

Of course, Theorems~\ref{Thm1} and~\ref{Thm2} can be combined
to disprove the existence of eigenvalues in a union of 
the \emph{unbounded} regions 
of the complex plane initially covered by the theorems separately.
From this perspective, our theorems are an improvement 
upon Cuenin \cite{Cuenin_2017}, 
who disproves the existence of eigenvalues
in \emph{compact} regions only.
Moreover, our results are much more explicit and quantitative.
On the other hand, our method seems to be restricted 
to the three-dimensional situation,
while the results of~\cite{Cuenin_2017} are stated in all dimensions
greater than or equal to two. 

\begin{Remark}
We note that~$V$ obeys~\eqref{hypothesis1} (respectively, \eqref{hypothesis2})
if, and only if, the adjoint~$V^*$ does. 
Noticing additionally that $f(\bar{\lambda},m)=f(\lambda,m)$, 
we see that any $\lambda \in \Com$ 
satisfying~\eqref{hypothesis1} or~\eqref{hypothesis2}
is not in the residual spectrum of~$H_V$ either.  	
Consequently, $\lambda$~is either in the continuous spectrum
or in the resolvent set of~$H_V$.
Finally, let us notice that 
$
  \sigma_\mathrm{ess}(H_V)
  = \sigma_\mathrm{ess}(H_0)
  = (-\infty,-m]\cup[+m,+\infty)
$
under the conditions stated in Theorems~\ref{Thm1} and~\ref{Thm2}.
\end{Remark}

Our approach to establish Theorems~\ref{Thm1} and~\ref{Thm2}
is based on two ingredients.
The main technique is the Birman-Schwinger principle 
stating that if~$\lambda$ is in the point spectrum of
the differential operator~$H_V$, 
then~$-1$ is in the spectrum of the (formal) integral operator
\begin{equation}\label{BS}
  W^{1/2} \, (H_{0}-\lambda)^{-1} \, U \, W^{1/2} 
  \,,
\end{equation}
where $V=UW$ is the polar decomposition of~$V$
(\ie~$W$ is the absolute value of~$V$ and~$U$ is unitary).
In this paper we make an effort to justify the principle
under the minimal regularity assumption that~$V$ is
in a sense relatively form-bounded with respect to~$H_0$.
Moreover, we rigorously cover the embedded eigenvalues
$\lambda \in (-\infty,-m]\cup[+m,+\infty)$.
Once the Birman-Schwinger principle is established, 
it is enough to show that the norm of~\eqref{BS} is less than one 
in order to establish Theorems~\ref{Thm1} and~\ref{Thm2}.
To this aim we use the second ingredient based on
a careful estimate of the resolvent integral kernel of~$H_0$.
Putting these two together, we end up with
the left-hand sides of the conditions~\eqref{hypothesis1} and~\eqref{hypothesis2}
being precisely upper bounds to the norm of~\eqref{BS}
obtained additionally with help of Sobolev-type inequalities.

The same approach has been recently employed in \cite[Sec.~2]{FKV}
to disprove the existence of eigenvalues of Schr\"odinger operators
with form-subordinated complex potentials (\cf~\cite[Thm.~1]{FKV}).
In the present case of Dirac operator, however, 
the method does not seem to lead to uniform results in~$\lambda$.

The rest of the paper consists of a precise definition of~$H_V$
as a pseudo-Friedrichs extension in Section~\ref{Sec.def},
the analysis of the free resolvent in Section~\ref{Sec.free},
the justification of the Birman-Schwinger principle in Section~\ref{Sec.BS}
and its application to the proofs of Theorems~\ref{Thm1} and~\ref{Thm2}
in Section~\ref{Sec.proof}.

%----------------------------------------%
\section{The pseudo-Friedrichs extension}\label{Sec.def}
%----------------------------------------%
%
We consider~$H_{V}$ as a perturbation of~$H_{0}$,
the latter being the operator that acts as~\eqref{Dirac} 
with $V=0$ and has the domain
$$
  \Dom(H_{0}) := 
  \left\{
  \psi \in \Hilbert :
  \nabla\psi \in \Hilbert^3
  \right\}
  .
$$
It is well known that~$H_0$ is self-adjoint
and that $C_0^\infty(\Real^3;\Com^4)$ is a core of~$H_0$.

Notice that $H_0^2 = (-\Delta+m^2) I_{\Com^4}$, 
where $-\Delta+m^2$ is the self-adjoint Schr\"odinger operator
in $\sii(\Real^3)$ with the usual domain $H^2(\Real^3)$.
The absolute value of~$H_0$ thus equals 
$(H_0^2)^{1/2} = \sqrt{-\Delta+m^2} I_{\Com^4}$,
which is again a self-adjoint operator 
when considered on the domain $H^1(\Real^3;\Com^4)$.
The form domain of $\sqrt{-\Delta+m^2}$ equals 
the fractional Sobolev space $H^{1/2}(\Real^3)$,
\cf~\cite[Sec.~7.11]{LL}.
Notice that $C_0^\infty(\Real^3)$ is dense in $H^{1/2}(\Real^3)$,
\cf~\cite[Sec.~7.14]{LL}. 
Clearly, $\sqrt{-\Delta+m^2} \geq \sqrt{-\Delta}$. 

We always assume that $|V| \in L^2_\mathrm{loc}(\Real^3)$
and that~$V$ is relatively form-bounded 
with respect to the massless~$H_0$ in the following sense:
There exist numbers $a \in (0,1)$ and $b\in\Real$ such that,
for all $\psi \in C_0^\infty(\Real^3)$,
\begin{equation}\label{Ass}
  \int_{\Real^3} |V(x)| |\psi(x)|^2 \, \der x
  \leq a \int_{\Real^3} |\sqrt[4]{-\Delta}\,\psi(x)|^2 \, \der x
  + b \int_{\Real^3} |\psi(x)|^2 \, \der x
  \,.
\end{equation}
We recall that $|V(x)|$ denotes the operator norm 
of the matrix~$V(x)$ in~$\Com^4$.
Similarly, $|\psi(x)|$ denotes the norm of the vector $\psi(x) \in \Com^4$. 
The double norm~$\|\cdot\|_\mathcal{X}$ 
is reserved for functional spaces~$\mathcal{X}$.
For instance, $\|\psi\|_{\Hilbert} = \||\psi|\|_{\sii(\Real)}$.
The symbol $(\cdot,\cdot)_\mathcal{X}$ stands for an inner product
in a Hilbert space~$\mathcal{X}$.

By the pseudo-Friedrichs extension \cite[Thm.~VI.3.11]{Kato}%
\footnote{Kato assumes $a<1/2$ for an abstract non-symmetric operator~$V$,
but it is straightforward to check that his proof works with $a<1$
(as in the symmetric case) in our special case 
when~$V$ is a multiplication by matrix and~\eqref{Ass} is assumed.}
(see also \cite{Veselic_2008} for more recent developments),
there exists a unique closed extension~$H_V$ of the operator sum $H_0+V$,
where~$V$ is understood as the multiplication operator 
by the matrix~$V$ in~$\Hilbert$
with initial domain $C_0^\infty(\Real^3;\Com^4)$. 
More specifically, the operator~$H_V$ satisfies
\begin{equation}\label{representation}
  (\phi,H_V\psi)_{\Hilbert}
  = \big(G_0^{1/2}\phi,H_0 G_0^{-1} G_0^{1/2}\psi\big)_{\Hilbert}
  + \int_{\Real^3} \phi(x)^* V(x) \psi(x) \, \der x
\end{equation}
for all $\psi \in \Dom(H_V) \subset H^{1/2}(\Real^3;\Com^4)$
and $\phi \in H^{1/2}(\Real^3;\Com^4)$.	
Here $G_0 := (H_0^2)^{1/2} + (a^{-1} b + \delta) I_{\Com^4}$ 
with any positive~$\delta$ 
(which does not influence the definition of~$H_V$). 
Notice that~$H_0$ and~$G_0$ commute
(in the usual sense for unbounded operators),
so the first term on the right-hand side of~\eqref{representation}
equals $(\phi,H_0\psi)_{\Hilbert}$ if $\psi \in H^1(\Real^3;\Com^4)$
Consequently, $H_V$~realises in a sense the form sum of~$H_0$ and~$V$. 
We have $\Dom({H_V}^{*}) \subset H^{1/2}(\Real^3;\Com^4)$
and $i\eta \not\in \sigma(H_V)$ for all real~$\eta$
with sufficiently large~$|\eta|$.

\begin{Remark}
By Kato's inequality
$
  \sqrt{-\Delta} \geq (2/\pi) |x|^{-1}
$ 
valid in the sense of quadratic forms in $\sii(\Real^3)$,
\cf~\cite[Rem.~V.5.12]{Kato},
we see that the potentials~$V$ satisfying the pointwise inequality
$$ 
  |V(x)| \leq a \, \frac{2}{\pi}\frac{1}{|x|} + b
$$ 
for almost every $x \in \Real^3$ obey the hypothesis~\eqref{Ass}. 
\end{Remark}

Another sufficient condition is given by the following proposition
(more optimal results can be obtained by using the Lorentz spaces).

\begin{Proposition}\label{Prop.sufficient}
Let $|V| = v_1 + v_2$ with $v_1 \in L^3(\Real^3)$ and $v_2 \in L^\infty(\Real^3)$
and assume that 
$$
  \|v_1\|_{L^3	(\Real^3)} < (2\pi^2)^{1/3}
  \,.
$$
Then~\eqref{Ass} holds true.  
\end{Proposition}
\begin{proof}
For every $\psi \in C_0^\infty(\Real^3)$, we have
$$
  \big\|v_1^{1/2}\psi\big\|_{\sii(\Real^3)}^2
  \leq \|v_1\|_{L^3(\Real^3)} 
  \, \big\|\psi\big\|_{L^3(\Real^3)}^{2}
  \leq (2\pi^2)^{-1/3} 
  \, \|v_1\|_{L^3(\Real^3)}
  \, \big\|\sqrt[4]{-\Delta}\,\psi\big\|_{L^2(\Real^3)}^{2}
  \,,
$$
where the first estimate follows by H\"older inequality
and the second inequality quantifies the Sobolev-type embedding 
$\dot{H}^{1/2}(\Real^3) \hookrightarrow L^3(\Real^3)$,
\cf~\cite[Thm.~8.4]{LL}.
At the same time, 
$
  \|v_2^{1/2}\psi\|_{\sii(\Real^3)}^2 
  \leq b \, \|\psi\|_{\sii(\Real^3)}^2
$
with some non-negative constant~$b$.
\end{proof}
\begin{Remark}
Since~\eqref{Ass} allows~$V$ to be merely defined through a quadratic form,
one can in principle introduce~$H_V$ under the weaker hypothesis 
$|V| \in L^1_\mathrm{loc}(\Real^3)$.
\end{Remark}
%

%---------------------------%
\section{The free resolvent}\label{Sec.free}
%---------------------------%
%
Using the well known trick 
$$
  (H_0-z)^{-1} = (H_0+z) (H_0^2-z^2)^{-1}
$$
and the knowledge of the integral kernel of 
the free Schr\"odinger operator $-\Delta + m^2$,
the integral kernel of the free resolvent $(H_{0}-z)^{-1}$
can be written down explicitly: 
\begin{equation}\label{Green}
  (H_0-z)^{-1}(x,x')
  = \frac{e^{-\sqrt{m^2-z^2}\;\!|x-x'|}}{4\pi\,|x-x'|}
  \left(
  \frac{i \alpha\cdot(x-x')}{|x-x'|^2}
  + \sqrt{m^2-z^2} \ \frac{i \alpha\cdot(x-x')}{|x-x'|}
  + m \alpha_4 + z
  \right)
\end{equation}
for every $z \not\in \sigma(H_0)$.
Here and in the sequel we choose the principal branch of the square root.
From now on, we also usually suppress writing 
the identity operators in the formulae;
therefore, we simply write~$z$ instead of $z I_\Hilbert$
(respectively, $z I_{\Com^4}	$)
on the left-hand (respectively, right-hand)
side of~\eqref{Green} and elsewhere. 

Given any matrix $M \in \Com^{4 \times 4}$,
we use the notation $|M|_\mathrm{HS}:=\sqrt{\tr(M^*M)}$
for the Hilbert-Schmidt (or Frobenius) norm of~$M$.
Recall that $|M| \leq |M|_\mathrm{HS}$.

\begin{Lemma}[Hilbert-Schmidt norm]\label{Lem.HS}
For almost every $x,x' \in \Real^3$, one has
\begin{multline*}
  \big|(H_0-z)^{-1}(x,x')\big|_\mathrm{HS}^2
  \\
  = 4 \, \frac{e^{-2\,\Re\sqrt{m^2-z^2}\,|x-x'|}}{(4\pi)^2 |x-x'|^4}
  \left(1+2\, \Re\sqrt{m^2-z^2} \, |x-x'|
  +2 \left[ \big(\Re\sqrt{m^2-z^2}\big)^2 + (\Re z)^2\right] |x-x'|^2\right)
  .
\end{multline*}
\end{Lemma}
\begin{proof}
Writing
$$  
  \big|(H_0-z)^{-1}(x,x')\big|_\mathrm{HS}^2
  = \frac{e^{-2\,\Re\sqrt{m^2-z^2}\,|x-x'|}}{(4\pi)^2 |x-x'|^6}
  \, \tr(A^*A)
  \,,
$$
where
\begin{equation}\label{A}
  A := a \cdot\alpha + a_4 \alpha_4 + a_0 I_{\Com^4}
\end{equation}
with
$$
  a := i \left(1+\sqrt{m^2-z^2}\,|x-x'|\right) (x-x')
  , 
  \qquad
  a_4 := m \, |x-x'|^2
  \,,
  \qquad
  a_0 := z \, |x-x'|^2
  \,,
$$
the proof reduces to straightforward manipulations with the Dirac matrices.
Using the anticommutation relations~\eqref{anti}, 
we immediately arrive at
$$
  A^* A 
  = (|a|^2 + a_4^2 + |a_0|^2) I_{\Com^4}
  + B 
  \,,
$$	
where 
$$
  B := b \cdot\alpha
  + b_4 \, \alpha_4 
  + c \cdot i\alpha\alpha_4
$$
with
$$
\begin{aligned}
  b &:= 2 \Re(\overline{a} a_0) (x-x')
  = 2 \left[
  \Im z - \big(\Im z \, \Re\sqrt{m^2-z^2} + \Re z \ \Im\sqrt{m^2-z^2} \big) |x-x'|
  \right] 
  |x-x'|^2 (x-x')
  \,,
  \\
  b_4 &:= 2 a_4 \Re(a_0) 
  = 2 m \, \Re z \, |x-x'|^4
  \,, 
  \\
  c &:= -2 \Im(a) a_4 (x-x')
  = -2 m \big( 1 + \Re\sqrt{m^2-z^2} \, |x-x'| \big) |x-x'|^2 (x-x')
  \,.
\end{aligned}
$$
Since $\tr(\alpha_\mu)=0$ for $\mu \in \{1,\dots,4\}$
as well as $\tr(\alpha_4\alpha_k)=0$ for $k \in \{1,\dots,3\}$
and $\tr(I_{\Com^4})=4$, we see that also the matrix~$B$ is traceless.
Consequently, 
using in addition the explicit expressions for $a_1$, $a_2$ and $a_3$,
we obtain
$$
\begin{aligned}
  \frac{1}{4}\tr(A^*A) 
  &= (|z|^2+m^2) |x-x'|^4 
  + \left|1+\sqrt{m^2-z^2}\,|x-x'|\right|^2 |x-x'|^2
  \\
  &= \left(|z|^2+m^2+|m^2-z^2|\right) |x-x'|^4 
  + 
  + 2 \Re \sqrt{m^2-z^2}\,|x-x'|^3 
  + |x-x'|^2
  \,.
\end{aligned}
$$
It remains to use the identity
$
  |z|^2+m^2+|m^2-z^2|
  = 2 \left[\big(\Re\sqrt{m^2-z^2}\big)^2+(\Re z)^2\right]
$.
\end{proof}

The present paper extensively uses the following explicit bound
on the Hilbert-Schmidt norm of the free resolvent.

\begin{Lemma}\label{Lem.bound}
For almost every $x,x' \in \Real^3$, one has
\begin{equation*}
  \big|(H_0-z)^{-1}(x,x')\big|_\mathrm{HS}^2
  \leq c_1^2 \, \frac{1}{|x-x'|^4} 
  + c_2^2 \, (\Re z)^2 \, \frac{e^{-2\,\Re\sqrt{m^2-z^2}\,|x-x'|}}{|x-x'|^2} 
  \,, 
\end{equation*}
where
$$
  c_1 := \frac{\sqrt{1+e^{-1}+2e^{-2}}}{2\pi} 
  \,, \qquad
  c_2 := \frac{\sqrt{2}}{2\pi}
  \,.
$$
\end{Lemma}
\begin{proof}
The estimate follows from Lemma~\ref{Lem.HS}
by using the elementary bounds 
$r e^{-r} \leq e^{-1}$ and $r^2 e^{-r} \leq 4 e^{-2}$
valid for every $r \geq 0$.
\end{proof}
%

%---------------------------------------%
\section{The Birman-Schwinger principle}\label{Sec.BS}
%---------------------------------------%
%
For almost every $x \in \Real^3$, 
let us introduce the matrix $W(x) := (V(x)^* V(x))^{1/2}$,
the absolute value of~$V(x)$.
Using the polar decomposition for matrices, 
we have 
$$
  V(x) = U(x) W(x) = U(x) W(x)^{1/2} W(x)^{1/2}
  \,,
$$
where~$U(x)$ is a unitary matrix on~$\Com^4$. 
Notice that $|V(x)| = |W(x)| = |W(x)^{1/2}|^2$,
where we recall that~$|\cdot|$ stands for the operator norm of a matrix.

Given any $z\not\in\sigma(H_0)$,
we introduce the \emph{Birman-Schwinger operator}
\begin{equation}%\label{BS}
  K_z :=
  \big[W^{1/2}G_0^{-1/2}\big]
  \big[G_0(H_{0}-z)^{-1}\big] 
  \big[(UW^{1/2})^*G_0^{-1/2}\big]^*   
  \,,  
\end{equation}
where~$G_0$ is the shifted absolute value of~$H_0$ introduced in~\eqref{Ass}.
With an abuse of notation, we denote by the same symbols~$W^{1/2}$ and~$U$ 
the maximal multiplication operators in~$\Hilbert$ 
generated by the matrix-valued function $x \mapsto W(x)$ 
and $x \mapsto U(x)$, respectively.
It follows from~\eqref{Ass} that, for every $\psi \in \Hilbert$,
$$
  \|W^{1/2}\psi\|_\Hilbert^2 \leq a \, \|G_0^{1/2}\psi\|_\Hilbert^2 
  \qquad\mbox{and}\qquad
  \|(UW^{1/2})^*\psi\|_\Hilbert^2 \leq a \, \|G_0^{1/2}\psi\|_\Hilbert^2 
  \,.
$$
Consequently, $\|W^{1/2}G_0^{-1/2}\|_{\Hilbert\to\Hilbert}^2 \leq a$ 
and $\|(UW^{1/2})^*G_0^{-1/2}\|_{\Hilbert\to\Hilbert}^2 \leq a$. 
Hence, $K_z$~is a well defined bounded operator on~$\Hilbert$ 
(as a composition of three bounded operators)
and one has the rough bound
\begin{equation}\label{rough}
  \|K_z\|_{\Hilbert\to\Hilbert} 
  \leq a \, \|G_0(H_{0}-z)^{-1}\|_{\Hilbert\to\Hilbert}
  \leq a \, \sup_{\xi\in(-\infty,-m)\cup(+m,\infty)}
  \frac{|\xi|+a^{-1}b+\delta}{|\xi-z|}
  \,.
\end{equation}
Using the commutativity of~$H_0$ and~$G_0$,
the operator~$K_z$ admits the familiar form~\eqref{BS}  
provided that~$V$ is bounded.
But we insist working under the minimal regularity assumption~\eqref{Ass}
in this section.
 
By $\varphi \in L_0^2(\Real^3;\Com^4)$ in the lemma below
we mean $\varphi \in L^2(\Real^3;\Com^4)$  
and that $\supp \varphi$ is compact.

\begin{Lemma}[Birman-Schwinger principle]\label{Lem.BS}
Assume~\eqref{Ass}. 
Let $H_V\psi=\lambda\psi$ with some $\lambda \in \Com$ and $\psi\in\Dom(H_V)$, 
and set $\phi := W^{1/2} \psi$.
One has $\phi \in \Hilbert$ and $\phi\not=0$ 
if $\psi\not=0$.
\begin{enumerate}
\item[\emph{(i)}]
If $\lambda \not\in \sigma(H_0)$, then $K_{\lambda} \phi = - \phi$.
\item[\emph{(ii)}]
If $\lambda \in \sigma(H_0)$, then
%
%\begin{equation}\label{BS}
$
\displaystyle
  \lim_{\eps \to 0^\pm} (\varphi,K_{\lambda + i\eps} \phi)_\Hilbert 
  = - (\varphi,\phi)_\Hilbert
$
for every $\varphi \in L_0^2(\Real^3;\Com^4)$.
%\end{equation}
%
\end{enumerate}
\end{Lemma}
\begin{proof}
By~\eqref{Ass}, $\phi \in \Hilbert$ if 
$\psi \in \Dom(H_V) \subset H^{1/2}(\Real^3;\Com^4)$.
From~\eqref{representation} we deduce that $H_V\psi=H_0\psi$
provided that $\phi=0$.
But then~$\lambda$ is an eigenvalue of~$H_0$ (which is impossible) unless $\psi=0$.  

If $\lambda \not\in \sigma(H_0)$, then
$$
\begin{aligned}
  (\varphi,K_\lambda\phi)_\Hilbert 
  &= \big([W^{1/2}G_0^{-1/2}]^*\varphi,
  [G_0(H_{0}-\lambda)^{-1}] [(UW^{1/2})^*G_0^{-1/2}]^*\phi\big)_\Hilbert
  \\
  &= \big( [UW^{1/2})^*G_0^{-1/2}]
  [G_0(H_{0}-\lambda)^{-1}]^* [W^{1/2}G_0^{-1/2}]^*  \varphi,
  W^{1/2}\psi\big)_\Hilbert
  \\
  &= \int_{\Real^3} \eta(x)^* V(x) \psi(x) \, \der x  
\end{aligned}
$$
for every $\varphi \in \Hilbert$ with
$
  \eta := G_0^{-1/2}[G_0(H_{0}-\lambda)^{-1}]^* [W^{1/2}G_0^{-1/2}]^*\varphi
  \in H^{1/2}(\Real^3;\Com^4)
  \,.
$
Using~\eqref{representation}, it follows that
$$
\begin{aligned}
  (\varphi,K_\lambda\phi)_\Hilbert 
  &= (\eta,H_V\psi)_\Hilbert	 
  - \big(G_0^{1/2}\eta,H_0 G_0^{-1} G_0^{1/2}\psi\big)_{\Hilbert}
  \\
  &= \lambda (\eta,\psi)_\Hilbert
  - \big(G_0^{1/2}\eta,H_0 G_0^{-1} G_0^{1/2}\psi\big)_{\Hilbert}
  \\
  &= \lambda \big(G_0^{1/2}\eta,G_0^{-1}G_0^{1/2}\psi\big)_\Hilbert
  - \big(G_0^{1/2}\eta,H_0 G_0^{-1} G_0^{1/2}\psi\big)_{\Hilbert}
  \\
  &= - \big(G_0^{1/2}\eta,(H_0-\lambda) G_0^{-1} G_0^{1/2}\psi\big)_{\Hilbert}
  \\
  &= - \big([W^{1/2}G_0^{-1/2}]^*\varphi,
  G_0(H_{0}-\lambda)^{-1}(H_0-\lambda) G_0^{-1} G_0^{1/2}\psi\big)_{\Hilbert}
  \\
  &= - \big([W^{1/2}G_0^{-1/2}]^*\varphi,
  G_0^{1/2}\psi\big)_{\Hilbert}
  \\ 
  &= - \big(\varphi,W^{1/2}\psi\big)_{\Hilbert}
  = - (\varphi,\phi)_{\Hilbert}
\end{aligned}
$$
for every $\varphi \in \Hilbert$.
This proves that $K_\lambda\phi=-\phi$ and therefore~(i).

If $\lambda \in \sigma(H_0)$, then there exists $\eps_0 > 0$
such that $\lambda+\eps \not\in \sigma(H_0)$
for all real~$\eps$ satisfying $0 < |\eps| < \eps_0$.
Now let us assume that $\varphi \in L_0^2(\Real^3;\Com^4)$.
As above, we have 
$$
\begin{aligned}
  (\varphi,K_{\lambda+i\eps}\phi)_\Hilbert 
  &= - \big([W^{1/2}G_0^{-1/2}]^*\varphi,
  G_0(H_{0}-\lambda-i\eps)^{-1}(H_0-\lambda) G_0^{-1} G_0^{1/2}\psi\big)_{\Hilbert}
  \\
  &= - (\varphi,\phi)_{\Hilbert}
  - i\eps \, \big([W^{1/2}G_0^{-1/2}]^*\varphi,
  G_0(H_{0}-\lambda-i\eps)^{-1} G_0^{-1} G_0^{1/2}\psi\big)_{\Hilbert}
  \\
  &= - (\varphi,\phi)_{\Hilbert}
  - i\eps \, \big(\varphi,
  W^{1/2}(H_{0}-\lambda-i\eps)^{-1}\psi\big)_{\Hilbert}
\end{aligned}
$$
for every $\varphi \in \Hilbert$.
Let us show that the last inner product vanishes as $\eps \to 0$.
Using Lemma~\ref{Lem.bound}, we have 
$$
  \big|\big(\varphi,
  W^{1/2}(H_{0}-\lambda-i\eps)^{-1}\psi\big)_{\Hilbert}\big|
  \leq  
  c_1 \, I_1 + c_2 \, |\lambda| \, I_2(\eps)
  \,,
$$
where
$$
\begin{aligned}
  I_1 &:= \int_{\Real^3\times\Real^3} 
  \frac{|\varphi(x)|\,|W(x)|^{1/2}\,|\psi(x')|}{|x-x'|^2} \, \der x \, \der x'
  \,,
  \\
  I_2(\eps) &:= \int_{\Real^3\times\Real^3} 
  |\varphi(x)|\,|W(x)|^{1/2} \, 
  \frac{e^{-\,\Re\sqrt{m^2-(\lambda+i\eps)^2}\,|x-x'|}}{|x-x'|} 
  \, |\psi(x')|
  \, \der x \, \der x'
  \,.
\end{aligned}
$$
The first integral is estimated as follows:
$$
\begin{aligned}
  I_1 
  &\leq (2\pi)^{2/3} \||\varphi||W|^{1/2}\|_{L^1(\Real^3)} 
  \||\psi|\|_{L^3(\Real^3)}
  \\
  &\leq (2\pi)^{2/3} \||\varphi||W|^{1/2}\|_{L^1(\Real^3)} 
  \, (2\pi^2)^{-1/6}
  \big\|\sqrt[4]{-\Delta}\,|\psi|\big\|_{L^2(\Real^3)}
  \\
  &\leq (2\pi)^{2/3} \||\varphi|\|_{L^2(\Real^3)}  
  \, \|\chi_\Omega|V|\|_{L^1(\Real^3)}^{1/2}
  \, (2\pi^2)^{-1/6}
  \big\|\sqrt[4]{-\Delta}\,|\psi|\big\|_{L^2(\Real^3)}
\end{aligned}
$$
with $\Omega := \supp\varphi$.
Here the first bound is due to 
the Hardy-Littlewood-Sobolev inequality \cite[Thm.~4.3]{LL},
the second estimate is 
a Sobolev-type inequality \cite[Thm.~8.4]{LL} quantifying 
the embedding $\dot{H}^{1/2}(\Real^3) \hookrightarrow L^3(\Real^3)$
and the last bound is the Schwarz inequality.
Notice that the integrals on the last time are all finite,
in particular because $\psi \in H^{1/2}(\Real^3;\Com^4)$
and $|V| \in L_\mathrm{loc}^1(\Real^3)$. 
To estimate the second integral, 
we proceed as in~\cite[proof of Lem.~2]{FKV}:
$$
  I_2(\eps) \leq \| |\varphi| \|_{L^2(\Real^3)}
  \|M_\eps\|_{L^2(\Real^3) \to L^2(\Real^3)}
  \| |\psi| \|_{L^2(\Real^3)} 
  \,,
$$
where $M_\eps$~is the integral operator with kernel
$$
  M_\eps(x,x') := 
  \chi_{\Omega}(x) \, |W(x)|^{1/2} \, 
  \frac{e^{-\,\Re\sqrt{m^2-(\lambda+i\eps)^2}\,|x-x'|}}{|x-x'|} 
  \,.
$$
Using the Hilbert-Schmidt norm, we have
$$
\begin{aligned}
  \|M_\eps\|_{L^2(\Real^3) \to L^2(\Real^3)}^2
  \leq  \|M_\eps\|_\mathrm{HS}^2
  &= \int_{\Omega\times\Real^3} 
  |V(x)| \, 
  \frac{e^{-\,2\Re\sqrt{m^2-(\lambda+i\eps)^2}\,|x-x'|}}{|x-x'|^2} 
  \, \der x \, \der x'
  \\
  &= \frac{2\pi}{\Re\sqrt{m^2-(\lambda+i\eps)^2}}
  \int_\Omega |V(x)| \, \der x
  \,,
\end{aligned}
$$ 
where the integral is again finite because of the local integrability of~$V$.
It remains to realise that
$$
  \Re\sqrt{m^2-(\lambda+i\eps)^2} \sim
  \begin{cases}
    |\eps|^{1/2} & \mbox{if} \quad \lambda^2=m^2 \ \& \ m\not=0 \,,
    \\
    |\eps| & \mbox{otherwise}. 
  \end{cases}
$$
Hence $\eps I_2(\eps) = O(|\eps|^{1/2})$ as $\eps \to 0$.
Putting the estimates together, we have proved
$$
\begin{aligned}
  \frac{|(\varphi,K_{\lambda+i\eps}\phi+\phi)_\Hilbert|}{\|\varphi\|_{\Hilbert}}
  &\leq \eps \, \|\chi_\Omega|V|\|_{L^1(\Real^3)}^{1/2} \,
  \||\psi|\|_{H^{1/2}(\Real^3)}
  \left(
  c_1 2^{1/2} \pi^{1/3}  
  + \frac{c_2 |\lambda| (2\pi)^{1/2}}
  {\big(\Re\sqrt{m^2-(\lambda+i\eps)^2}\big)^{1/2}} 
  \right)
  %\\
  \xrightarrow[\eps \to 0]{} 0
\end{aligned}
$$
for every given $\varphi \in L_0^2(\Real^3;\Com^4)$.
This shows~(ii) and concludes the proof of the lemma.
\end{proof}

The preceding lemma is a precise statement of one side
of the Birman-Schwinger principle 
under the minimal regularity assumption~\eqref{Ass}.
It says that if~$\lambda$ is an eigenvalue of~$H_V$,
then $-1$~is an eigenvalue of an integral equation related to~$K_\lambda$. 
If $\lambda \not\in\sigma(H_0)$ the converse implication also holds,
but it is not generally true if $\lambda \in\sigma(H_0)$,
\cf~\cite[Sec.~III.2]{SiQF}, 
and it is not needed for the purpose of this paper.
In fact, we exclusively use the following corollary of Lemma~\ref{Lem.BS}.

\begin{Corollary}\label{Corol.BS}.
Assume~\eqref{Ass} and let $\lambda \in \sigma_\mathrm{p}(H_V)$. 
\begin{enumerate}
\item[\emph{(i)}]
If $\lambda \not\in \sigma(H_0)$, 
then $\|K_{\lambda}\| \geq 1$.
\item[\emph{(ii)}]
If $\lambda \in \sigma(H_0)$, then
$
\displaystyle
  \liminf_{\eps\to 0^\pm} \|K_{\lambda+i\eps}\|_{\Hilbert\to\Hilbert} \geq 1
$.
\end{enumerate}
\end{Corollary}
\begin{proof}
Let $\lambda \in \sigma_\mathrm{p}(H_V)$, 
let~$\psi$ be a corresponding eigenfunction
and set $\phi = W^{1/2}\psi \not= 0$.

If $\lambda \not\in \sigma(H_0)$, 
then the statement~(i) of Lemma~\ref{Lem.BS} implies
$$
  \|\phi\|_\Hilbert^2 \, \|K_\lambda\|_{\Hilbert\to\Hilbert} 
  \geq |(\phi,K_\lambda\phi)_\Hilbert| = \|\phi\|_\Hilbert^2
  \,,
$$
from which the claim~(i) immediately follows.

If $\lambda \in \sigma(H_0)$, 
we set $\phi_n := \xi_n \phi$ for every positive~$n$,
where~$\xi_n(x):=\xi(x/n)$ and $\xi \in C_0^\infty(\Real^3)$
is a usual cut-off function satisfying
$\xi(x)=1$ for $|x| \leq 1$ and $\xi(x)=0$ for $|x| \geq 2$.
As above, we write
$$
  \|\phi_n\|_\Hilbert \|\phi\|_\Hilbert \, 
  \|K_{\lambda+i\eps}\|_{\Hilbert\to\Hilbert} 
  \geq |(\phi_n,K_{\lambda+i\eps}\phi)_\Hilbert| 
  \,,
$$
Taking the limit $\eps \to 0^\pm$,
the statement~(ii) of Lemma~\ref{Lem.BS} implies
$$
  \|\phi_n\|_\Hilbert \|\phi\|_\Hilbert \, 
  \liminf_{\eps\to 0^\pm} \|K_{\lambda+i\eps}\|_{\Hilbert\to\Hilbert} 
  \geq |(\phi_n,\phi)_\Hilbert|
  \,.
$$
The desired claim~(ii) then follows by taking the limit $n\to\infty$.
\end{proof}

To disprove the existence of eigenvalues in a complex region,
it is thus enough to show that the norm of the Birman-Schwinger
operator is strictly less than one there.

%---------------%
\section{Proofs}\label{Sec.proof}
%---------------%
%
To implement the above idea, we first strengthen our hypotheses about~$V$
and establish a more precise estimate 
on the norm of~$K_z$ as compared to~\eqref{rough}.

\begin{Lemma}\label{Lem1}
Let $|V| \in L^3(\Real^3)$.
For every $z \not\in \sigma(H_0)$, 
$$
  \|K_z\|_{\Hilbert\to\Hilbert}
  \leq \left(\frac{\pi}{2}\right)^{1/3} \sqrt{1+e^{-1}+e^{-2}} \, 
  \sqrt{ 
  1 + \frac{(\Re z)^2}{\big(\Re\sqrt{m^2-z^2}\big)^2} 
  }
  \
  \||V|\|_{L^3(\Real^3)} 
  \,.
$$
\end{Lemma}
\begin{proof}
By the elementary bound $r^2 e^{-r} \leq 4 e^{-2}$ valid for every $r \geq 0$,
we get from Lemma~\ref{Lem.bound} the estimate
\begin{equation*}
  \big|(H_0-z)^{-1}(x,x')\big|_\mathrm{HS}^2
  \leq \left( 
  c_1^2 + \tilde{c}_2^2 \frac{(\Re z)^2}{\big	(\Re\sqrt{m^2-z^2}\big)^2} 
  \right) 
  \frac{1}{|x-x'|^4} 
\end{equation*}
for almost every $x,x' \in \Real^3$
with $\tilde{c}_2 := c_2 e^{-2} = \sqrt{2e^{-2}}/(2\pi)$.
Consequently, for every $\phi,\psi \in \Hilbert$,
$$
\begin{aligned}
  |(\phi,K_z\psi)_{\Hilbert}| 
  &\leq 
  \sqrt{ 
  c_1^2 + \tilde{c}_2^2 \frac{(\Re z)^2}{\big(\Re\sqrt{m^2-z^2}\big)^2} 
  } 
  \int_{\Real^3 \times \Real^3} 
  \frac{|V(x)|^{1/2} \, |\phi(x)| \, |V(x')|^{1/2} \, |\psi(x')| }
  {|x-x'|^2}
  \, \der x \, \der x'
  \\	
  &\leq 
  \sqrt{ 
  c_1^2 + \tilde{c}_2^2 \frac{(\Re z)^2}{\big(\Re\sqrt{m^2-z^2}\big)^2} 
  } 
  \ 
  2^{2/3}\pi^{4/3} \, 
  \||V|\|_{L^3(\Real^3)}  
  \, 
  \|\phi\|_\Hilbert \|\psi\|_\Hilbert
  \,,
\end{aligned}
$$
where the second estimate is due to the Hardy-Littlewood-Sobolev
inequality with optimal constants \cite[Thm.~4.3]{LL}
together with the H\"older inequality.
We eventually use $\tilde{c}_2 \leq c_1$, 
just to make the final result look more elegant.
\end{proof}

Further strengthening hypothesis about~$V$,
we can also get a bound which is independent of the mass~$m$
and has a different behaviour in~$z$.
\begin{Lemma}\label{Lem2}
Let $|V| \in L^3(\Real^3) \cap L^{3/2}(\Real^3)$.
For every $z \not\in \sigma(H_0)$, 
$$
  \|K_z\|_{\Hilbert\to\Hilbert}
  \leq \left(\frac{\pi}{2}\right)^{1/3} \sqrt{1+e^{-1}+e^{-2}} \, 
  \||V|\|_{L^3(\Real^3)} 
  + \frac{2^{17/6}}{3\pi^{2/3}} \, |\Re z| \,
  \||V|\|_{L^{3/2}(\Real^{3})} 
  \,.
$$
\end{Lemma}
\begin{proof}
For every $\phi,\psi \in \Hilbert$,
we have 
$$
\begin{aligned}
  |(\phi,K_z\psi)_{\Hilbert}| 
  &\leq 
  \int_{\Real^3 \times \Real^3} |V(x)|^{1/2} \, |\phi(x)|
  \left(\frac{c_1}{|x-x'|^2} + \frac{c_2 \, |\Re z|}{|x-x'|} \right)
  |V(x')|^{1/2} \, |\psi(x')| 
  \, \der x \, \der x'
  \\
  &\leq 
  \left(
  c_1 \, 2^{2/3}\pi^{4/3} \, 
  \||V|\|_{L^3(\Real^3)}  
  + c_2 \, |\Re z| \, \frac{2^{10/3}}{3 \pi^{1/3}} \, 
  \||V|\|_{L^{3/2}(\Real^{3})} 
  \right)
  \|\phi\|_\Hilbert \|\psi\|_\Hilbert
  \,.
\end{aligned}
$$
Here the first inequality follows by Lemma~\ref{Lem.bound}
(where the exponential is estimated by~$1$)
and the second estimate is due to the Hardy-Littlewood-Sobolev
inequality with optimal constants \cite[Thm.~4.3]{LL}
together with the H\"older inequality.
\end{proof}

Theorem~\ref{Thm1} (respectively, Theorem~\ref{Thm2})
follows as a consequence of Corollary~\ref{Corol.BS}
and Lemma~\ref{Lem1} (respectively, Lemma~\ref{Lem2}).

%---------------------------%
\subsection*{Acknowledgment}
%---------------------------%
%
The research of D.K.\ was partially supported 
by the GACR grant No.\ 18-08835S
and by FCT (Portugal) through project PTDC/MAT-CAL/\-4334/\-2014.
 

%\newpage
%\vfill
%--------------%
% BIBLIOGRAPHY %
%--------------%
%
%\addcontentsline{toc}{section}{References}
\bibliography{bib}
\bibliographystyle{amsplain}

\end{document}